\documentclass[a4paper]{article}
\usepackage{amsfonts}
\usepackage{amssymb}
\usepackage{amsmath}
\usepackage{latexsym}
\usepackage{graphicx}
\usepackage{calrsfs}
\newtheorem{theorem}{Theorem}[section]
\newtheorem{corollary}[theorem]{Corollary}

\newtheorem{example}[theorem]{Example}

\newtheorem{remark}[theorem]{Remark}

\begin{document}

\parbox{1mm}

\begin{center}
{\bf {\sc \Large Toeplitz Operators on Poly-analytic Spaces via
Time-Scale Analysis}}
\end{center}

\vskip 12pt

\begin{center}
{\bf Ondrej HUTN\'IK and M\'aria HUTN\'IKOV\'A}\footnote{{\it
Mathematics Subject Classification (2010):} Primary 47B35,
Secondary 47G30, 47L80, 42C40,
\newline {\it Key words and phrases:} Time-scale analysis, poly-analytic function,
Toeplitz operator, operator symbol, Bergman space}
\end{center}


\hspace{5mm}\parbox[t]{12cm}{\fontsize{9pt}{0.1in}\selectfont\noindent{\bf
Abstract.} 
This is a review paper based on the series of our papers devoted
to a structure of true-poly-analytic Bergman function spaces over
the upper half-plane in the complex plane and to a detailed study
of properties of Toeplitz operators with separate symbols acting
on them via methods of time-scale analysis.} \vskip 24pt

\section{Introduction}

Analytic functions are the main object of classical complex
analysis. One definition of analytic functions is in terms of
Cauchy-Riemann operator
$\partial_{\overline{z}}:=\frac{1}{2}(\frac{\partial}{\partial
x}+\mathrm{i}\frac{\partial}{\partial y})$ with $z=x+\mathrm{i}y$.
Then $f$ is analytic in some (simply connected bounded or
unbounded) domain $\Omega$ of the complex plane $\mathbb{C}$ iff
$\partial_{\overline{z}} f=0$ on $\Omega$. A natural extension of
this definition is to iterate the Cauchy-Riemann operator which
yields a notion of poly-analytic function, i.e., $f$ is
poly-analytic (or, analytic of order $n$) in $\Omega$ iff
$\partial_{\overline{z}}^n f=0$ on $\Omega$. In this paper we will
consider $\Omega=\Pi$ -- the upper half-plane in the complex
plane. Poly-analytic functions of order $n$ which are not
poly-analytic of any other lower order are called
true-poly-analytic of order $n$. For further reading on
poly-analytic functions we refer to~\cite{balk}.


It was recently observed in~\cite{abreu} that \textit{the
true-poly-analytic Bergman space of order $k+1$} (over the upper
half-plane) \textit{may be alternatively viewed as the space of
wavelet transforms with Laguerre functions of order $k$}. This
perhaps unexpected result enables us to see poly-analytic spaces
and related operators of complex analysis from a new perspective
of time-scale analysis. In this paper we are especially interested
in Toeplitz operators on (true)-poly-analytic Bergman spaces. A
recent interest of this topic may be seen e.g. in~\cite{CL}.

The paper is a survey of results from~\cite{hutnik}-\cite{HH},
where further results, details, and proofs can be found. The main
ingredients of the whole theory (on the hand of time-scale
analysis) are the affine group $\mathbb{G}$ of
orientation-preserving linear transformations of the real line,
and a parameterized family of admissible affine coherent states
(wavelets) whose on Fourier transform side are defined via
Laguerre functions. In this context we first describe the
structure of the space of wavelet transforms of Hardy-space
functions inside the Hilbert space $L_2(\mathbb{G},
\mathrm{d}\nu_L)$, and link this construction with the
intertwining property of induced representation of $\mathbb{G}$.
Indeed, this study provides a time-scale approach to poly-analytic
spaces as explained in~\cite{abreu}. Then the obtained results are
applied to the study of behavior of Toeplitz operators on wavelet
subspaces (i.e., true-poly-analytic Bergman spaces on $\Pi$) which
provides an interesting generalization of the classical Toeplitz
operators acting on (weighted) Bergman spaces studied
in~\cite{vasilevskibook}. Thus, we deal with boundedness and other
properties of Toeplitz operators and their algebras with symbols
as individual coordinates of the underlying affine group.

\section{Time-scale approach to poly-analytic spaces: a description}\label{secaffine}

\paragraph{Affine group and its induced representations.}
The \textit{affine group} $\mathbb{G}$ consists of all
transformations $A_{u,v}$ of the real line $\mathbb{R}$ of the
type $A_{u,v}(x):= ux+v$, $x\in\mathbb{R}$, where $u>0$,
$v\in\mathbb{R}$. Indeed, writing
$$\mathbb{G}=\{\zeta=(u,v); \,\,u>0, v\in\mathbb{R}\},$$ one has the
multiplication law on $\mathbb{G}$ of the form $\zeta_1\diamond
\zeta_2 = (u_1,v_1)\diamond(u_2,v_2) = (u_1u_2,u_1v_2+v_1)$. With
respect to the multiplication $\diamond$ the group $\mathbb{G}$ is
non-commutative with the identity element $e=(1,0)$, and locally
compact Lie group on which the left-invariant Haar measure is
given by $\mathrm{d}\nu_L(\zeta)=u^{-2}\,\mathrm{d}u\mathrm{d}v$.
The usual identification of the group $\mathbb{G}$ with the upper
half-plane $\Pi=\{\zeta=v+\mathrm{i}u; \,v\in\mathbb{R}, u>0\}$ in
the complex plane $\mathbb{C}$ equipped with the hyperbolic metric
and the corresponding (hyperbolic) measure $\mathrm{d}\nu_L$ will
also be used. Then $L_2(\mathbb{G},\mathrm{d}\nu_L)$ denotes the
Hilbert space of all square-integrable complex-valued functions on
$\mathbb{G}$ with respect to the measure $\mathrm{d}\nu_L$.

The affine group $\mathbb{G}$ may be decomposed as a semi-direct
product $\mathbb{G}=\textsf{N}\ltimes \textsf{A}$, where
$\textsf{N}=\{(1,v); v\in\mathbb{R}\}$ is the abelian normal
closed subgroup, and the quotient group $\textsf{A}$ is isomorphic
to the one-parameter closed subgroup $\{(u,0); u>0\}\cong
\mathbb{R}_+$. Thus, if $H$ is a closed subgroup of $\mathbb{G}$
and $X=\mathbb{G}/H$ is the corresponding left-homogeneous space,
we may induce representations of $\mathbb{G}$ in the subspaces
which depend on $X=\mathbb{G}/H$ with $H=\{e\}$, $H=\textsf{A}$
and $H=\textsf{N}$, respectively. Indeed,

\begin{itemize}
\item[(i)] $X=\mathbb{G}/\{e\}=\mathbb{G}$ -- a character of the
subgroup $\{e\}$ induces a left-regular representation of
$\mathbb{G}$ on $L_2(X) = L_2(\mathbb{G}, \mathrm{d}\nu_L)$ in the
form $$[\Lambda(u,v)F](x,y) :=
F\left(\frac{x}{u},\frac{y-v}{u}\right);$$ \item[(ii)]
$X=\mathbb{G}/\textsf{N}\cong \textsf{A}=\mathbb{R}_+$ -- a
character of the subgroup $\textsf{N}$ induces a co-adjoint
representation of $\mathbb{G}$ on
$L_2(\mathbb{G}/\textsf{N})=L_2(\mathbb{R}_+)$ in the form
$$\left[\rho\left(u,v\right)f\right]\left(x\right)
:= \mathrm{e}^{-2\pi \mathrm{i}\frac{v}{x}}
f\left(\frac{x}{u}\right);$$ \item[(iii)]
$X=\mathbb{G}/\textsf{A}\cong \textsf{N}=\mathbb{R}$ -- a
character of the subgroup $\textsf{A}$ induces a quasi-regular
representation of $\mathbb{G}$ on
$L_2(\mathbb{G}/\textsf{A})=L_2(\mathbb{R})$ in the form $$
\left[\pi\left(u,v\right)f\right]\left(y\right) :=
\frac{1}{\sqrt{u}}\,f\left(\frac{y-v}{u}\right);$$
\end{itemize} whenever $(u,v), (x,y)\in\mathbb{G}$.

The Hilbert space $L_{2}\left( \mathbb{R}\right)$ under the action
$\pi$ contains precisely two closed proper invariant subspaces
$H_{2}\left(\mathbb{R}\right)$ and
$H_{2}\left(\mathbb{R}\right)^{\bot}$, called the Hardy and
conjugate Hardy spaces, respectively, such that $L_2(\mathbb{R}) =
H_2(\mathbb{R})\oplus H_2(\mathbb{R})^\bot$. Thus, $\pi$ is a
reducible representation on $L_{2}\left( \mathbb{R}\right)$, and
we can decompose it into two irreducible representations, such
that $ \pi(u,v) = \pi^+(u,v) \oplus\pi^-(u,v)$. From it follows
that only the Hardy space $H_2(\mathbb{R})$ is considered,
although the discussion and further results are equally valid for
the conjugate Hardy space $ H_{2}\left( \mathbb{R}\right)^{\bot}$.
From the action on signals (we identify a signal with an element
$f \in L_2(\mathbb{R})$) we observe that $\mathbb{G}$ consists
precisely of the transformations we apply to a signal:
\textit{translation} (time-shift) by an amount $v$, and
\textit{zooming in or out} by the factor $u$. Hence, the group
$\mathbb{G}$ naturally relates to the geometry of signals.

There is an intertwining operator between the co-adjoint
representation $\rho$ and quasi-regular representation $\pi$. This
is the Fourier transform $\mathcal{F}: L_{2}(\mathbb{R})\to
L_{2}(\mathbb{R})$ in the form $\mathcal{F}\{f\}(\xi) :=
\hat{f}(\xi)=\int_{\mathbb{R}} f(x)\,\mathrm{e}^{-2\pi
\mathrm{i}x\xi}\,\mathrm{d}x$, since it uses the characters which
induce those representations. Also there exists an intertwining
operator between $\pi$ and $\Lambda$ given by the identity
$[W_\psi f](u,v) := \langle f, \pi(u,v) \psi\rangle$, $f\in
L_2(\mathbb{R})$, which provides a starting point for time-scale
analysis on $\mathbb{R}$. Usually it is desirable to make this map
unitary as well, and this is expressed by the following resolution
of identity
$$\langle f,g\rangle = \int_\mathbb{G} \langle f,
\pi(\zeta)\psi\rangle \langle \pi(\zeta)\psi,
g\rangle\,\mathrm{d}\nu_L(\zeta),$$ also known as the Calder\'on
reproducing formula. To achieve this the mother wavelet $\psi\in
L_2(\mathbb{R})$ shall be \textit{admissible}: for the affine
group this is equivalent to
$$
\int_{\mathbb{R}_+} |\hat{\psi}(u)|^2\,\frac{\mathrm{d}u}{u} = 1.
$$

\paragraph{Wavelets from Laguerre functions and their subspaces.}
Here we describe an alternative approach to true-poly-analytic
Bergman spaces using wave-lets built from Laguerre functions.
Thus, for $k\in\mathbb{Z}_+$ consider the parameterized family of
admissible wavelets $\psi^{(k)}$ on $\mathbb{R}$ defined on the
Fourier transform side as follows
$$\hat{\psi}^{(k)}(\xi)=\chi_+(\xi)\sqrt{2\xi}\,\ell_k(2\xi),$$
where $\ell_n(x):=\mathrm{e}^{-x/2} L_n(x)$ are the (simple)
Laguerre functions with $L_n(x)$ being the Laguerre polynomial of
order $n\in\mathbb{Z}_+$ and $\chi_+$ the characteristic function
of the positive half-line. Then according to the Calder\'on
reproducing formula
\begin{align*}
f(v) & =
\int_{\mathbb{R}_+}\Bigl(\bigl(D_u\psi^{(k)}\bigr)*\bigl(D_u\psi^{(k)}\bigr)*
f\Bigr)(v)\,\frac{\mathrm{d}u}{u^2}
\end{align*}for all $f\in
H_2(\mathbb{R})$, where $*$ is the usual convolution on
$L_2(\mathbb{R})$. For each $k\in\mathbb{Z}_+$ define the
subspaces $A^{(k)}$ of $L_2(\mathbb{G},\mathrm{d}\nu_L)$ as
follows
\begin{align*}
A^{(k)} & := \left\{ [W_k
f](u,v)=\Bigl(f*\bigl(D_u\psi^{(k)}\bigr)\Bigr)(v); \,f\in
H_2(\mathbb{R})\right\}.
\end{align*}Indeed, $W_k f$ are exactly the
continuous wavelet transforms of functions $f\in H_2(\mathbb{R})$
with respect to wavelets $\psi^{(k)}$. Consequently, $A^{(k)}$
will be referred to as \textit{wavelet subspaces} of
$L_2(\mathbb{G}, \mathrm{d}\nu_L)$. Note that it is possible to
consider the "conjugate" wavelet
$\hat{\bar{\psi}}^{(k)}(\xi)=\hat{\psi}^{(k)}(-\xi)$, the
"conjugate" wavelet subspaces
\begin{align*}
\bar{A}^{(k)} & := \left\{ [W_{\bar{k}}
f](u,v)=\Bigl(f*\bigl(D_u\bar{\psi}^{(k)}\bigr)\Bigr)(v); \,f\in
H_2(\mathbb{R})^\bot\right\}
\end{align*} and to build up the
theory in this setting. In what follows we will state the results
only for $A^{(k)}$, similar results may be stated for its
"conjugate" counterpart $\bar{A}^{(k)}$.

\begin{remark}\rm
Note that poly-analytic Bergman spaces and introduced wavelet
subspaces share intriguing patterns that may prove usable. A
deeper study of this connection is given in the recent
paper~\cite{abreu}: the important and interesting observation of
that paper is that for $k\in\mathbb{N}$ \textit{the spaces
$A^{(k-1)}$ of continuous wavelet transforms of Hardy space
functions with respect to wavelets from Laguerre functions
coincide with the true-poly-analytic Bergman spaces of order $k$
on the upper half-plane} (symmetrically, $\bar{A}^{(k-1)}$
corresponds to the space of all true-poly-anti-analytic functions
of order $k$ from $L_2(\Pi)$). This allows to study these objects
of complex analysis using techniques of time-scale analysis.
\end{remark}

\begin{figure}
\begin{center}
\includegraphics[scale=0.9]{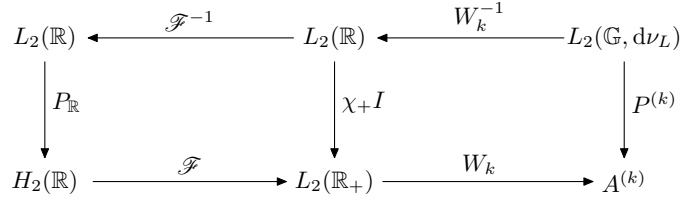}
\end{center}\caption{Relationship among the introduced spaces and operators}\label{fig_operatorsI}
\end{figure}

The relationship among the introduced spaces $A^{(k)}$ of wavelet
transforms of $H_2(\mathbb{R})$-functions, and the unitary
operators of continuous wavelet transform $W_k$ and the Fourier
transform $\mathcal{F}$ is schematically described on
Figure~\ref{fig_operatorsI}. For each $k\in\mathbb{Z}_+$ the
spaces $A^{(k)}$ are the reproducing kernel Hilbert spaces.
Explicit formulas for their reproducing kernels
$$K_\zeta^{(k)}(\eta)=\left\langle \pi(\eta)\psi^{(k)},
\pi(\zeta)\psi^{(k)}\right\rangle$$ and orthogonal projections
$P^{(k)}: L_2(\mathbb{G}, \mathrm{d}\nu_L)\to A^{(k)}$ are
described in~\cite{hutnik}.


\paragraph{Structural results.}
In accordance with the representation of $L_2(\mathbb{G},
\mathrm{d}\nu_L)$ as tensor product in the form
$$L_2(\mathbb{G}, \mathrm{d}\nu_L(\zeta)) = L_{2}(\mathbb{R}_+, u^{-2}\mathrm{d}u)
\otimes L_{2}(\mathbb{R}, \mathrm{d}v)$$ with $\zeta = (u,v) \in
\mathbb{G}$, we consider the unitary operator $$U_{1}=
(I\otimes\mathcal{F}): L_{2}(\mathbb{R}, u^{-2}\mathrm{d}u)
\otimes L_{2}(\mathbb{R}, \mathrm{d}v) \longrightarrow
L_{2}(\mathbb{R}_+, u^{-2}\mathrm{d}u) \otimes L_{2}(\mathbb{R},
\mathrm{d}\omega).$$ For the purpose to "linearize" the hyperbolic
measure $\mathrm{d}\nu_L$ onto the usual Lebes\-gue plane measure
we introduce the unitary operator
$$U_{2}: L_{2}(\mathbb{R}_+,\,u^{-2}\mathrm{d}u)\otimes
L_{2}(\mathbb{R},\,\mathrm{d}\omega) \longrightarrow
L_{2}(\mathbb{R}_+,\,\mathrm{d}x)\otimes
L_{2}(\mathbb{R},\,\mathrm{d}y)$$ given by
$$U_{2}: F(u,\omega) \longmapsto
\frac{\sqrt{2|y|}}{x} F\left(\frac{x}{2|y|},y\right).$$
Immediately we get the following theorem describing the structure
of $A^{(k)}$ inside $L_2(\mathbb{G},\mathrm{d}\nu_L)$.

\begin{theorem}[\cite{hutnik}, Theorem~2.1]\label{thm1}
The unitary operator $U=U_{2}U_{1}$ gives an isometrical
isomorphism of the space $L_{2}(\mathbb{G},\mathrm{d}\nu_L)$ onto
$L_{2}(\mathbb{R}_+,\,\mathrm{d}x)\otimes
L_{2}(\mathbb{R},\,\mathrm{d}y)$ under which
\begin{itemize}
\item[(i)] the wavelet subspace $A^{(k)}$ is mapped onto
$L_k\otimes L_2(\mathbb{R}_+)$, where $L_{k}$ is the rank-one
space generated by Laguerre function $\ell_{k}(x)$; \item[(ii)]
the orthogonal projection $P^{(k)}:
L_2(\mathbb{G},\mathrm{d}\nu_L)\to A^{(k)}$ is unitarily
equivalent to the following one $UP^{(k)}U^{-1} =
P_0^{(k)}\otimes\chi_+ I,$ where $P_0^{(k)}$ 
is the one-dimensional projection of $L_2(\mathbb{R}_+,
\mathrm{d}x)$ onto $L_k$.
\end{itemize}
\end{theorem}

\begin{remark}\rm
Let us mention that connection between certain spaces of wavelet
transforms and Bergman spaces is already well-known, see
e.g.~\cite{paul}: the Bergman transform $[B^\alpha F](u,v) =
u^{-\alpha-1/2} F(u,v)$ with $\alpha>0$ gives an isometrical
isomorphism of the space $L_2(\mathbb{G}, \mathrm{d}\nu_L)$ onto
$L_2(\mathbb{G}, u^{2\alpha-1} \mathrm{d}u\mathrm{d}v)$ under
which the space of continuous wavelet transforms of
$H_2(\mathbb{R})$-functions with respect to the Bergman wavelet
given by $\hat{\psi}^\alpha(\xi)=\chi_+(\xi) c_\alpha \xi^\alpha
\mathrm{e}^{-2\pi\xi}$ is mapped onto the weighted Bergman space
$\mathcal{A}_{2\alpha-1}(\mathbb{G})$. Here, $c_\alpha$ is a
certain normalization constant.
\end{remark}

However, we may say more about the connection between the wavelet
subspaces and Hardy spaces which reads as follows. Let us mention
that the orthogonal projection $P_{\mathbb{R}}$ of
$L_2(\mathbb{R})$ onto $H_2(\mathbb{R})$ is called the
\textit{Szeg\"{o} projection}.

\begin{theorem}\label{thm2}
The unitary operator
$V=(I\otimes\mathcal{F}^{-1})U_{2}(I\otimes\mathcal{F})$ gives an
isometrical isomorphism of the space
$L_{2}(\mathbb{G},\mathrm{d}\nu_L)$ onto
$L_{2}(\mathbb{R}_+,\,\mathrm{d}x)\otimes
L_{2}(\mathbb{R},\,\mathrm{d}y)$ under which
\begin{itemize}
\item[(i)] $A^{(k)}$ and $H_2(\mathbb{R})$ are connected by
the formula $V\left(A^{(k)}\right) = L_k\otimes H_2(\mathbb{R});$  
\item[(ii)] $P^{(k)}$ and $P_\mathbb{R}$ are connected by the
formula $VP^{(k)}V^{-1} = P_0^{(k)}\otimes P_\mathbb{R}$.
\end{itemize}
\end{theorem}

\begin{figure}
\begin{center}
\includegraphics[scale=0.9]{operators2.1}
\end{center}\caption{Visualizing the results of Theorem~\ref{thm1}
and Theorem~\ref{thm2}}\label{fig_operatorsII}
\end{figure}

The diagram on Figure~\ref{fig_operatorsII} schematically
describes all the relations among the constructed operators and
spaces appearing in the above two theorems. A detailed analysis of
the construction and origin of unitary operators describing the
structure of wavelet subspaces from the viewpoint of induced
representations of $\mathbb{G}$ is done in~\cite{EH}. It was shown
that these unitary maps have the following properties related to
group representations:
\begin{enumerate}
\item[(i)] they intertwine respective representations of the
affine group $\mathbb{G}$; \item[(ii)] they provide a spatial
separation of the irreducible components of the affine group's
representations.
\end{enumerate}
Indeed, these properties make the unitary maps useful for
characterization of $A^{(k)}$ inside the space $L_2(\mathbb{G},
\mathrm{d}\nu_L)$.

\begin{remark}\rm
The suggested time-scale (or, more general time-frequency) point
of view was recently successfully used in another Abreu's
paper~\cite{abreu} to obtain a complete characterization of all
lattice sampling and interpolating sequences in the
Segal-Bargmann-Fock space of poly-analytic functions or,
equivalently, of all lattice vector-valued Gabor frames and
vector-valued Gabor Riesz sequences with Hermite functions for
$L_2(\mathbb{R},\mathbb{C}^n)$. This again underlines a new, and
perhaps unexpected, connection between poly-analytic functions and
time-frequency analysis having a great potential in various
applications.
\end{remark}

\section{Toeplitz operators on poly-analytic spaces}

Toeplitz operators form one of the most significant classes of
concrete operators because of their importance both in pure and
applied mathematics and in many other sciences. In the context of
true-poly-analytic Bergman spaces (or, equivalently, wavelet
subspaces) for a given bounded function $a$ on $\mathbb{G}$ define
the \textit{Toeplitz operator} $T_a^{(k)}: A^{(k)}\to A^{(k)}$
with symbol $a$ usually as $T_{a}^{(k)} := P^{(k)}M_{a}$, where
$M_a$ is the operator of pointwise multiplication by $a$ on
$L_2(\mathbb{G},\mathrm{d}\nu_L)$ and $P^{(k)}$ is the orthogonal
projection from $L_2(\mathbb{G},\mathrm{d}\nu_L)$ onto $A^{(k)}$.
In fact, this provides a mapping between the same wavelet
subspaces. It is worth noting that in the case of many wavelet
subspaces (parameterized by $k$) other Toeplitz- and Hankel-type
operators may be defined, e.g.
\begin{align*}
T_a^{(k,l)} & := P^{(k)}M_a P^{(l)}, \\
h_a^{(k,l)} & := \bar{P}^{(k)}M_a P^{(l)}, \\
H_a^{(k,l)} & := \left(I-\sum_{j=0}^{k}P^{(j)}\right)M_a P^{(l)}.
\end{align*}In what follows we restrict our attention only to the case
$k=l$. The mapping $a\mapsto T_a^{(k)}$ is then interpreted as the
quantization rule "on the level $k$".

\subsection{Unitarily equivalent images of Toeplitz operators for
symbols depending on $\Im\zeta$}

It was observed in several cases~\cite{vasilevskibook} that
Toeplitz operators can be transformed into pseudo-differential
operators by means of certain unitary maps constructed as an exact
analog of the Bargmann transform mapping the Segal-Bargmann-Fock
space $F_2(\mathbb{C}^n)$ of Gaussian square-integrable entire
functions on complex $n$-space onto $L_2(\mathbb{R}^n)$. Via this
mapping Toeplitz operator $T_a^{(k)}: A^{(k)}\to A^{(k)}$ can be
identified with certain pseudo-differential operator
$\mathfrak{T}_a^{(k)}: L_2(\mathbb{R})\to L_2(\mathbb{R})$ which
provides an analog of the Berezin reducing of Toeplitz operators
with anti-Wick symbols on the Fock space $F_2(\mathbb{C}^n)$ to
Weyl pseudo-differential operators on $L_2(\mathbb{R}^n)$. This
construction may be better seen if we apply this procedure to
operator symbols $a(u,v)$ that are only depending on individual
variables. Indeed, for the case when $a = a(u)$ depends only on
the first spatial variable of $\mathbb{G}$ the operator
$T_a^{(k)}$ is simply a multiplication operator with explicitly
computable symbol.

\begin{theorem}[\cite{hutnikIII}, Theorem~3.2]\label{CTO1}
Let $(u,v)\in \mathbb{G}$. If a measurable symbol $a=a(u)$ does
not depend on $v$, then $T_{a}^{(k)}$ acting on $A^{(k)}$ is
unitarily equivalent to the multiplication operator
$\mathfrak{A}_a^{(k)} = \gamma_{a,k} I$ acting on
$L_2(\mathbb{R}_+)$, where the function $\gamma_{a,k}:
\mathbb{R}_+\to\mathbb{C}$ is given by
\begin{equation}\label{gamma1}
\gamma_{a,k}(\xi) =
\int_{\mathbb{R}_{+}}a\left(\frac{u}{2\xi}\right)\ell_k^2(u)\,\mathrm{d}u,
\quad \xi\in\mathbb{R}_+.
\end{equation}
\end{theorem}

\begin{example}\rm\label{exm12}
Given a point $\lambda_0\in\mathbb{R}_+$ we have
$$\gamma_{\chi_{[0,\lambda_0]},k}(\xi) =
\chi_+(\xi)\int_{\mathbb{R}_+}
\chi_{[0,\lambda_0]}\left(\frac{u}{2\xi}\right)\ell_k^2(u)\,\mathrm{d}u
= \chi_+(\xi)\int_{0}^{2\lambda_0 \xi} \ell_k^2(u)\,\mathrm{d}u.$$
Immediately, for $\lambda_0=\frac{1}{2}$ and $k=0$ we have
$\gamma_{\chi_{[0,1/2]},0}(\xi) = 1-\mathrm{e}^{-\xi}$,
$\xi\in\mathbb{R}_+$.
\end{example}

The function $\gamma_{a,k}$ is obtained by integrating a dilation
of a symbol $a=a(u)$ of $T_a^{(k)}$ against a Laguerre function of
order $k$. This result extends the result of Vasilevski for the
classical Toeplitz operators acting on the Bergman space (i.e.,
the case $k = 0$ in our notation) in very interesting way which
differs from the case of Toeplitz operators acting on weighted
Bergman spaces summarized in Vasilevski
book~\cite{vasilevskibook}. Moreover, the function $\gamma_{a,k}$
sheds a new light upon the investigation of main properties of the
corresponding Toeplitz operator $T_a^{(k)}$ with a symbol $a=a(u)$
(measurable and unbounded, in general), such as boundedness,
spectrum, invariant subspaces, norm value, etc. Furthermore, the
use unitarily equivalent images as model, or local representatives
permits us to study Toeplitz operators with much more general
symbols.

Since the function $\gamma_{a,k}: \mathbb{R}_+\to\mathbb{C}$ is
responsible for many interesting features of the corresponding
Toeplitz operator $T_a^{(k)}$, we present here certain interesting
and useful properties of $\gamma_{a,k}$ in what follows.

\begin{theorem}\label{thmder}
Let $(u,v)\in \mathbb{G}$. If $a=a(u)\in L_1(\mathbb{R}_+)\cup
L_\infty(\mathbb{R}_+)$ such that $\gamma_{a,k}(\xi)\in
L_\infty(\mathbb{R}_+)$, then for each $n=1,2,\dots$
$$\lim_{\xi\to+\infty}\frac{\mathrm{d}^n\gamma_{a,k}(\xi)}{\mathrm{d}\xi^n}
= 0$$ holds for each $k\in\mathbb{Z}_+$. Moreover, if $a=a(u)\in
C_b^\infty(\mathbb{R}_+)$ such that for each $n\in\mathbb{N}$
holds
$$\lim_{u\to+\infty} u^n\,a^{(n)}(u) = 0,$$ then also
$$\lim_{\xi\to 0} \xi^n \frac{\mathrm{d}^n \gamma_{a,k}(\xi)}{\mathrm{d}\xi^n} = 0$$ for each $n\in\mathbb{N}$ and each
$k\in\mathbb{Z}_+$.
\end{theorem}

The result of above theorem states that the behavior of
derivatives of $\gamma_{a,k}$ \textit{does not depend on
parameter} $k$. In fact, it depends only on behavior of the
corresponding symbol $a$ (or, its derivatives), but not on the
particularly chosen Laguerre functions. This is quite surprising
because, as we have already stated, the wavelet transforms with
Laguerre functions of order $k$ live, up to a multiplier
isomorphism, in the true-poly-analytic Bergman space of order $k$,
which is rather different from the classical Bergman space of
analytic functions. Thus, we have the remarkable observation that,
asymptotically, all the true-poly-analytic Bergman spaces have
"the same behavior". This result has some important consequences
in quantum physics, signal analysis and in the asymptotic theory
of random matrices, which are not yet completely understood.

\begin{remark}\rm\label{remSV}
The special case $n=1$ yields the following result: if $a\in
C_b^1(\mathbb{R}_+)$ with $\lim\limits_{u\to+\infty} u a'(u)=0$,
then for each $k\in \mathbb{Z}_+$ the function $\gamma_{a,k}$ is
\textit{slowly varying at infinity} (in the additive sense) and
\textit{slowly varying at zero} (in the multiplicative sense).
\end{remark}

Easily, for each $k\in\mathbb{Z}_+$ and each $a=a(u)\in
L_\infty(\mathbb{R}_+)$ we have
$$\sup_{\xi\in\mathbb{R}_+} |\gamma_{a,k}(\xi)| \leq \sup_{u\in\mathbb{R}_+} |a(u)| \int_{\mathbb{R}_+} \ell_k^2(u)\,\mathrm{d}u < +\infty,$$
i.e., $\gamma_{a,k}$ is bounded on $\mathbb{R}_+$ for each
$k\in\mathbb{Z}_+$. Moreover, in such a case of bounded symbol $a$
the function $\gamma_{a,k}(\xi)$ is also continuous in each finite
point $\xi\in\mathbb{R}_+$, and thus $\gamma_{a,k}\in
C_b(\mathbb{R}_+)$. However, as the following examples show
$\gamma_{a,k}$ may be bounded even for unbounded symbols.

\begin{example}\rm\label{exmunbounded}
(i) For unbounded symbol
$$a(u)=\frac{1}{\sqrt{u}}\sin\frac{1}{u}, \quad u\in\mathbb{R}_+,
$$ we have
$$\gamma_{a,1}(\xi) = \frac{\sqrt{2\pi}}{4}\,\mathrm{e}^{-2\sqrt{\xi}}
\left[\left(2\sqrt{\xi}-8\xi\right)\frac{\cos
2\sqrt{\xi}}{2\sqrt{\xi}}+\left(3-2\sqrt{\xi}\right)\frac{\sin
2\sqrt{\xi}}{2\sqrt{\xi}}\right] $$ for $\xi\in\mathbb{R}_+$,
which is a bounded function on $\mathbb{R}_+$. However, due to
computational limitations we can not say anything about the
boundedness of $\gamma_{a,k}(\xi)$ for arbitrary $k$.

(ii) For oscillating symbol $a(u)=\mathrm{e}^{2u \mathrm{i}}$ we
have again the bounded function
$$\gamma_{a,k}(\xi) =
\frac{(-1)^k}{(\xi-\mathrm{i})^{2k+1}}\sum_{j=0}^{k}
(-1)^{j}\left[{k\choose j}\right]^2 \xi^{2j+1}, \quad
\xi\in\overline{\mathbb{R}}_+.$$ Moreover, $\gamma_{a,k}(\xi)\in
C[0,+\infty]$ for each $k\in\mathbb{Z}_+$.
\end{example}

These examples motivate to study this interesting feature in more
detail considering unbounded symbols to have a sufficiently large
class of them common to all admissible $k$. For this purpose
denote by $L_{1}(\mathbb{R}_+, 0)$ the class of functions $a=a(u)$
such that $a(u)\,\mathrm{e}^{-\varepsilon u} \in
L_1(\mathbb{R}_+)$ for any $\varepsilon > 0$. 
For any such $L_{1}(\mathbb{R}_+, 0)$-symbol $a(u)$ define the
following averaging functions
\begin{align*}
C_{a}^{(1)}(u) = \int_0^u a(t)\,\mathrm{d}t, \quad C_{a}^{(m)}(u)
& = \int_0^u C_{a}^{(m-1)}(t)\,\mathrm{d}t, \quad m=2,3,\dots
\end{align*}The functions $C_a^{(m)}$ constitute a "sequence of iterated
integrals" of symbol $a$.

\begin{theorem}\label{thmmeans1}
Let $a=a(u)\in L_{1}(\mathbb{R}_+, 0)$.
\begin{itemize}
\item[(i)] If for any $m\in\mathbb{N}$ the function $C_a^{(m)}$
has the following asymptotic behavior
\begin{equation}\label{asymptoticC1}
C_{a}^{(m)}(u) = O(u^m),\,\,\,\, \textrm{as}\,\,\,\, u\to 0\,\,\,
\textrm{as well as}\,\,\,\, u\to +\infty,
\end{equation}
then for each $k\in\mathbb{Z}_+$ we have
$\sup_{\xi\in\mathbb{R}_+} |\gamma_{a,k}(\xi)|<+\infty$.

\item[(ii)] If for any $m,n\in\mathbb{N}$, any
$\lambda_1\in\mathbb{R}_+$ and any $\lambda_2\in (0,n+1)$ holds
\begin{equation}\label{asymptoticC4}
C_a^{(m)}(u) = O\left(u^{m+\lambda_1}\right),\,\,\,\,
\textrm{as}\,\,\,\, u\to 0, \end{equation} and
$$
C_a^{(n)}(u) = O\left(u^{n-\lambda_2}\right),\,\,\,\,
\textrm{as}\,\,\,\, u\to +\infty,$$ then for each
$k\in\mathbb{Z}_+$ we have $\lim\limits_{\xi\to+\infty}
\gamma_{a,k}(\xi) = 0 = \lim\limits_{\xi\to 0} \gamma_{a,k}(\xi)$.
\end{itemize}
\end{theorem}

\begin{remark}\rm
The condition~(\ref{asymptoticC1}) guarantees the boundedness of
the function $\gamma_{a,k}(\xi)$ at a neighborhood of
$\xi=+\infty$, as well as at a neighborhood of $\xi=0$. Observe
that if the both conditions in~(\ref{asymptoticC1}) hold for some
$m=m_0$, then they hold also for $m=m_0+1$. Indeed,
$$|C_a^{(m_0+1)}(u)| \leq \int_0^u |C_a^{(m_0)}(t)|\,\mathrm{d}t
\leq \textrm{const}\,\int_0^u t^{m_0}\,\mathrm{d}t \leq
\textrm{const}\,u^{m_0+1}.$$
\end{remark}

The main advantage of Theorem~\ref{thmmeans1} is that we need not
have an explicit form of the corresponding function $\gamma_{a,k}$
for an unbounded symbol $a=a(u)$ to decide about its boundedness.
Also, it gives the condition on the behavior of
$L_1(\mathbb{R}_+,0)$-\-sym\-bols such that the function
$\gamma_{a,k}(\xi)\in C[0,+\infty]$.

\begin{example}\rm\label{exmunbounded5}
For $\alpha>0$ and $\beta \in (0,1)$ consider the unbounded symbol
$$a(u)=u^{-\beta}\sin u^{-\alpha}, \quad u\in\mathbb{R}_+.$$ However, the function $a(u)$ is
continuous at $u=+\infty$ for all admissible values of parameters,
and therefore $\gamma_{a,k}(0) = a(+\infty) = 0$. On the other
side, it is difficult to verify the behavior of function
$\gamma_{a,k}(\xi)$ at the endpoint $+\infty$ by a direct
computation. Since
$$C_a^{(1)}(u) = \frac{u^{\alpha-\beta+1}}{\alpha} \cos u^{-\alpha}
+ O(u^{2\alpha-\beta+1}), \quad \textrm{as}\,\,\, u\to
0,$$
then for $\alpha>\beta$ the first condition
in~(\ref{asymptoticC4}) holds for $m=1$ and
$\lambda_1=\alpha-\beta$. By Theorem~\ref{thmmeans1} the function
$\gamma_{a,k}(\xi)$ is bounded.

\noindent If $\alpha\leq\beta$, then
$$C_a^{(m)}(u) =
O(u^{m\alpha-\beta+m}),\,\,\,\,\textrm{as}\,\,\,\, u\to 0.$$ Thus,
for each $\alpha\leq\beta$ there exists $m_0\in\mathbb{N}$ such
that $m_0\alpha>\beta$, and therefore the first condition
in~(\ref{asymptoticC4}) holds for $m=m_0$ and
$\lambda_1=m_0\alpha-\beta$, which guarantees that
$\gamma_{a,k}(\xi)$ is continuous at $\xi=0$. Thus, for all
parameters $\alpha>0$ and $\beta\in (0,1)$ the function
$\gamma_{a,k}(\xi)\in C[0,+\infty]$ for each $k\in\mathbb{Z}_+$.
\end{example}

We also mention an alternative way to the properties of
$\gamma_{a,k}$: using the explicit form of Laguerre polynomial we
may write
\begin{equation}\label{widetildegamma}
\gamma_{a,k}(\xi) = 2\xi \int_{\mathbb{R}_+}
a(u)\,\ell_k^2(2u\xi)\,\mathrm{d}u = \sum_{i=0}^k \sum_{j=0}^k
\kappa(k,i,j)\, \widetilde{\gamma}_{a,i+j}(\xi),
\end{equation}
where
$$\widetilde{\gamma}_{a,\lambda}(\xi) := (2\xi)^{\lambda+1} \int_{\mathbb{R}_+}
a(u)\,u^\lambda\,\mathrm{e}^{-2u\xi}\,\mathrm{d}u =
\xi^{\lambda+1} \int_{\mathbb{R}_+}
a(u/2)\,u^\lambda\,\mathrm{e}^{-u\xi}\,\mathrm{d}u,$$ and
$$\kappa(k,i,j) := \frac{(-1)^{i+j}}{i! j!} {k\choose i}
{k\choose j}.$$ It is immediate that the boundedness of each
function $\widetilde{\gamma}_{a,\lambda}$ for
$\lambda\in\{0,1,\dots,2k\}$ implies the boundedness of function
$\gamma_{a,k}$. Therefore, given $\lambda\in\mathbb{Z}_+$ and a
locally sumable function $a=a(u)$ we now introduce the weighted
means of symbol $a$ as follows
\begin{align*}
D_{a,\lambda}^{(1)}(u) = \int_0^u a(t/2) t^\lambda\,\mathrm{d}t,
\quad D_{a,\lambda}^{(m)}(u) & = \int_0^u
D_{a,\lambda}^{(m-1)}(t)\,\mathrm{d}t, \quad
m=2,3,\dots\end{align*}It is obvious that $D_{a,0}^{(m)}(u) = 2^m
C_a^{(m)}(u/2)$ for each $m\in\mathbb{N}$. In this setting the
result of Theorem~\ref{thmmeans1} reads as follows: if $a=a(u)\in
L_{1}(\mathbb{R}_+, 0)$ and for any $m\in\mathbb{N}$ the function
$D_{a,0}^{(m)}$ has the asymptotic behavior
$$D_{a,0}^{(m)}(u) = O(u^m),\,\,\,\, \textrm{as}\,\,\,\, u\to
0\,\,\, \textrm{as well as}\,\,\,\, u\to +\infty,$$ then each
function $\widetilde{\gamma}_{a,\lambda}$ is bounded on
$\mathbb{R}_+$ for every $\lambda\in\mathbb{Z}_+$, which means
that $\gamma_{a,k}$ is bounded on $\mathbb{R}_+$ for each
$k\in\mathbb{Z}_+$. Moreover, we may extend the above observation
about asymptotic behavior using any (positive) weight
$\lambda_0\in\mathbb{Z}_+$ appearing in weighted means
$D_{a,\lambda_0}^{(m)}$.

\begin{theorem}
Let $a=a(u)\in L_{1}(\mathbb{R}_+, 0)$. If for any $\lambda_0\in
\mathbb{Z}_+$ and for any $m\in\mathbb{N}$ the function
$D_{a,\lambda_0}^{(m)}$ has the following asymptotic behavior
$$D_{a,\lambda_0}^{(m)}(u) = O(u^{\lambda_0+m}),\,\,\,\,
\textrm{as}\,\,\,\, u\to 0\,\,\, \textrm{as well as}\,\,\,\, u\to
+\infty,$$ then $\gamma_{a,k}$ is bounded for each
$k\in\mathbb{Z}_+$.
\end{theorem}

Furthermore, if $a=a(u)$ is a bounded symbol on $\mathbb{G}$, then
the operator $T_a^{(k)}$ is clearly bounded on each $A^{(k)}$, and
for its operator norm holds $\|T_a^{(k)}\| \leq
\textrm{ess-sup}\,|a(u)|$. Thus, all spaces $A^{(k)}$,
$k\in\mathbb{Z}_+$, are naturally appropriate for Toeplitz
operators with bounded symbols. However, we may observe that the
result of Theorem~\ref{CTO1} suggests considering not only
$L_{\infty}(\mathbb{G},\mathrm{d}\nu_L)$-symbols, but also
\textit{unbounded} ones. In this case we obviously have

\begin{corollary}\label{corgamma}
The operator $T_a^{(k)}$ with a measurable symbol $a=a(u)$,
$u\in\mathbb{R}_+$, is bounded on $A^{(k)}$ if and only if the
corresponding function $\gamma_{a,k}(\xi)$ is bounded on
$\mathbb{R}_+$, and
$$\|T_a^{(k)}\| = \sup_{\xi \in\mathbb{R}_+}
|\gamma_{a,k}(\xi)|.$$
\end{corollary}

From this result we immediately have that \textit{all the obtained
results for boundedness of $\gamma_{a,k}$ in terms of iterated
integrals are, in fact, sufficient conditions for boundedness of
the corresponding Toeplitz operator $T_a^{(k)}$ on each
$A^{(k)}$.} The following result provides another criteria for
simultaneous boundedness of Toeplitz operators on each wavelet
subspaces.

\begin{theorem}
(i) Let $a=a(u)\in L_{1}(\mathbb{R}_+, 0)$ be non-negative almost
everywhere. If $T_a^{(0)}$ is bounded on $A^{(0)}$, then the
operator $T_a^{(k)}$ is bounded on $A^{(k)}$ for each
$k\in\mathbb{Z}_+$.

(ii) Let $C_a^{(m)}$ be non-negative almost everywhere for a
certain $m=m_0$. If $T_a^{(0)}$ is bounded on $A^{(0)}$, then the
operator $T_a^{(k)}$ is bounded on $A^{(k)}$ for each
$k\in\mathbb{Z}_+$.
\end{theorem}

Theorem states that under the assumption of non-negati\-vi\-ty of
a symbol $a\in L_1(\mathbb{R}_+, 0)$, or its mean $C_{a}^{(m)}$
for certain $m\in\mathbb{N}$, the boundedness of Toeplitz operator
$T_a = T_a^{(0)}$ on the Bergman space $\mathcal{A}_2(\Pi) =
A^{(0)}$ implies the boundedness of Toeplitz operator $T_a^{(k)}$
acting on $A^{(k)}$ for each $k\in\mathbb{Z}_+$. However, if $a\in
L_1(\mathbb{R}_+, 0)$, the question whether the boundedness of
$T_a^{(k_0)}$ on $A^{(k_0)}$ for certain $k_0\in \mathbb{N}$
implies the boundedness of $T_a^{(k)}$ acting on $A^{(k)}$ for
each $k\in\mathbb{Z}_+$ (smaller, or greater than $k_0$) is still
open. This is related to question whether the boundedness of
$T_a^{(k)}$ may happen only simultaneously for all
$k\in\mathbb{Z}_+$.

It is immediate that an unbounded symbol must have a sufficiently
sophisticated oscillating behavior at neighborhoods of the points
$0$ and $+\infty$ to generate a bounded Toeplitz operator. In what
follows we show that infinitely growing positive symbols cannot
generate bounded Toeplitz operators in general. For this purpose
for a non-negative function $a=a(u)$ put
$$\theta_{a}(u)=\inf_{t\in (0,u)} a(t)\,\,\,\,\textrm{and}\,\,\,\,\Theta_{a}(u)=\inf_{t\in (u/2,u)} a(t).$$

\begin{theorem}\label{thm3}
For a given non-negative symbol $a=a(u)$ if
\begin{equation}\label{limm_a,0}
\textrm{either}\,\,\,\lim_{u\to 0}
\theta_a(u)=+\infty,\,\,\,\textrm{or}\,\,\, \lim_{u\to +\infty}
\Theta_{a}(u)=+\infty,
\end{equation} then $T_a^{(k)}$
is unbounded on each $A^{(k)}$, $k\in\mathbb{Z}_+$.
\end{theorem}

\begin{example}\rm
For the family of non-negative symbols on $\mathbb{R}_+$ in the
form
$$a(u)=u^{-\beta}\ln^2 u^{-\alpha}, \quad \beta\in[0,1],\,
\alpha>0,$$ we have that for all admissible parameters holds
$\lim\limits_{u\to 0} \theta_a(u) = +\infty$, and thus $T_a^{(k)}$
is unbounded on $A^{(k)}$ for each $k\in\mathbb{Z}_+$.
\end{example}

In the following we provide an interesting example of symbols
$a,b$ for which $T_a^{(k)}, T_b^{(k)}$ are bounded, but
$T_{ab}^{(k)}$ is not on the whole scale of parameters $k$.

\begin{example}\rm
Let us consider two symbols on $\mathbb{R}_+$ in the form $$a(u) =
u^{-\beta} \sin u^{-\alpha}, \,\,\, \beta\in (0,1),\,\, \alpha
\geq \beta,\quad \textrm{and}\quad b(u) = u^{\tau} \sin
u^{-\alpha}, \,\,\, \tau\in (0,\beta).$$ Then $T_a^{(k)}$ is
bounded for each $k\in\mathbb{Z}_+$, and since $b(u)\in
C[0,+\infty]$, then $T_b^{(k)}$ is bounded for each
$k\in\mathbb{Z}_+$ as well. Put
$$c(u) = a(u)b(u) = \frac{u^{-\delta}}{2} - \frac{u^{-\delta}}{2}
\cos 2u^{-\alpha} = c_1(u)+c_2(u),$$ where $\delta=\beta-\tau \in
(0,1)$. Clearly, $c(u)$ is an unbounded symbol. However,
$T_{c_2}^{(k)}$ is bounded for each $k\in\mathbb{Z}_+$. Since
$$\theta_{c_1}(u) = \inf_{t\in (0,u)} \frac{1}{2 t^{\delta}} =
\frac{1}{2 u^{\delta}} \to +\infty,
\,\,\,\,\textrm{as}\,\,\,\,u\to 0,$$ then the operator
$T_{c_1}^{(k)}$ is unbounded for each $k\in\mathbb{Z}_+$. Thus,
the Toeplitz operator $T_{ab}^{(k)}$ is unbounded on $A^{(k)}$ for
each $k\in\mathbb{Z}_+$ showing that the semi-commutator
$\left[T_a^{(k)}, T_b^{(k)}\right)$ \textit{is not compact}.
\end{example}

Perhaps the most surprising feature of behavior of Toeplitz
operators on $A^{(k)}$ with symbols depending only on vertical
coordinate in the upper half-plane is appearance of certain
commutative algebras of Toeplitz operators on these spaces which
are practically unknown in the literature. Therefore, denote by
$L^{\{0,+\infty\}}_\infty(\mathbb{R}_+)$ the $C^*$-subalgebra of
$L_\infty(\mathbb{R}_+)$ which consists of all functions having
limits at the points $0$ and $+\infty$. For $k\in\mathbb{Z}_+$
denote by
$\mathcal{T}_k\left(L^{\{0,+\infty\}}_\infty(\mathbb{R}_+)\right)$
the $C^*$-algebra generated by all operators $T_a^{(k)}$ acting on
$A^{(k)}$ with symbols $a\in
L^{\{0,+\infty\}}_\infty(\mathbb{R}_+)$.

\begin{theorem}\label{thmlimits}
For $a\in L^{\{0,+\infty\}}_\infty(\mathbb{R}_+)$ the
corresponding functions $\gamma_{a,k}(\xi)$, $k\in\mathbb{Z}_+$,
possess the following properties
\begin{itemize}
\item[(i)] $\gamma_{a,k}(\xi)\in C[0,+\infty]$; \item[(ii)]
$\gamma_{a,k}(+\infty)=\lim\limits_{\xi\to +\infty}
\gamma_{a,k}(\xi) = \lim\limits_{u\to 0} a(u)=a(0)$; \item[(iii)]
$\gamma_{a,k}(0)=\lim\limits_{\xi\to 0} \gamma_{a,k}(\xi) =
\lim\limits_{u\to +\infty} a(u)=a(+\infty)$.
\end{itemize}
\end{theorem}

Indeed, the behavior of a bounded function $a(u)$ near the point
$0$, or $+\infty$ determines the behavior of function
$\gamma_{a,k}(\xi)$ near the point $+\infty$, or $0$,
respectively. An interesting observation is that the limits at
infinity and at zero of the function $\gamma_{a,k}$ \textit{are
completely independent of parameter $k$}, which is rather
surprising in this case and it again may be useful in various
contexts. The existence of limits of $a(u)$ at these endpoints
guarantees the continuity of $\gamma_{a,k}(\xi)$ on $[0,+\infty]$,
however this condition is not necessary even for bounded symbols
as the following example shows.

\begin{example}\rm
For $a(u)=\sin u$, $u\in\mathbb{R}_+$, we have
$$\gamma_{a,1}(\xi) = 2\xi \int_{\mathbb{R}_+} \sin u\,
\mathrm{e}^{-2u\xi}(1-2u\xi)^2\,\mathrm{d}u =
\frac{2\xi\,(1-16\xi^2+48\xi^4)}{(1+4\xi^2)^3}, \quad
\xi\in\mathbb{R}_+,$$ which yields
$$\lim_{\xi\to +\infty} \gamma_{a,1}(\xi) = \lim_{\xi\to 0}
\gamma_{a,1}(\xi) = 0.$$ On the other hand,
Example~\ref{exmunbounded}(i) provides an example of unbounded
symbol $a(u)$ such that the corresponding function $\gamma_{a,k}$
\textit{is continuous} on $[0,+\infty]$ for each
$k\in\mathbb{Z}_+$.
\end{example}

\begin{corollary}\label{coralgebra}
Each $C^*$-algebra
$\mathcal{T}_k\left(L^{\{0,+\infty\}}_\infty(\mathbb{R}_+)\right)$,
$k\in\mathbb{Z}_+$, is isometric and isomorphic to $C[0,+\infty]$.
The corresponding isomorphism is generated by the following
mapping $\tau^{(k)}: T_a^{(k)} \longmapsto \gamma_{a,k}(\xi)$.
\end{corollary}

According to this result the symbol algebra
$L^{\{0,+\infty\}}_\infty(\mathbb{R}_+)$ is an example of algebra
such that the operator algebra
$\mathcal{T}_k\left(L^{\{0,+\infty\}}_\infty(\mathbb{R}_+)\right)$
generated by Toeplitz operators $T_a^{(k)}$ with symbol $a\in
L^{\{0,+\infty\}}_\infty(\mathbb{R}_+)$ \textit{is commutative for
each} $k\in\mathbb{Z}_+$. An interesting question may be a
characterization of all such algebras $\mathcal{A}$ of symbols for
which the operator algebra $\mathcal{T}_k(\mathcal{A})$ is
commutative for each $k$. It seems to be a challenging problem.

Easily, continuity of function $\gamma_{a,k}$ on the whole
$\overline{\mathbb{R}}_+$ guarantees its boundedness, and
therefore the boundedness of the corresponding Toeplitz operator
$T_a^{(k)}$ on $A^{(k)}$. Moreover, $T_a^{(k)}$ with a symbol
$a=a(u)$ belongs to the algebra
$\mathcal{T}_k\left(L_\infty^{\{0,+\infty\}}(\mathbb{R}_+)\right)$
if and only if the corresponding function $\gamma_{a,k}(\xi)$
belongs to $C[0,+\infty]$. This means that the algebra
$\mathcal{T}_k\left(L_\infty^{\{0,+\infty\}}(\mathbb{R}_+)\right)$
also contains (bounded) Toeplitz operators whose (generally
unbounded) symbols $a(u)$ need not have limits at the endpoints
$0$ and $+\infty$. For instance, Example~\ref{exmunbounded}(ii)
provides oscillating symbols $a=a(u)$ for which Toeplitz operator
$T_a^{(k)}$ belongs to the algebra
$\mathcal{T}_k\left(L_\infty^{\{0,+\infty\}}(\mathbb{R}_+)\right)$
for the whole range of parameters $k$. Now we give an example of a
bounded oscillating symbol such that the bounded operator
$T_a^{(k)}$ does not belong to the algebra
$\mathcal{T}_k\left(L_\infty^{\{0,+\infty\}}(\mathbb{R}_+)\right)$.

\begin{example}\rm The function $a(u) = u^\mathrm{i} = \mathrm{e}^{\mathrm{i} \ln u}$,
$u\in\mathbb{R}_+,$ is oscillating near the endpoints $0$ and
$+\infty$, but it is bounded on $\mathbb{R}_+$, and therefore
$T_a^{(k)}$ is bounded for each $k\in\mathbb{R}_+$. Changing the
variable $x=2u\xi$ yields $$\gamma_{a,k}(\xi) =
2\xi\int_{\mathbb{R}_+} u^\mathrm{i} \ell_k^2(2u\xi)\,\mathrm{d}u
= (2\xi)^{-\mathrm{i}} \int_{\mathbb{R}_+} x^\mathrm{i}
\ell_k^2(x)\,\mathrm{d}x.$$ Since the last integral is a constant
depending on $k$, the function $\gamma_{a,k}(\xi)$ oscillates and
has no limit when $\xi\to 0$ as well as when $\xi\to +\infty$.
Thus, the bounded Toeplitz operator $T_a^{(k)}$ \textit{does not
belong} to the algebra
$\mathcal{T}_k\left(L_\infty^{\{0,+\infty\}}(\mathbb{R}_+)\right)$.
Hence not all oscillating symbols (even bounded and continuous)
generate an operator from
$\mathcal{T}_k\left(L_\infty^{\{0,+\infty\}}(\mathbb{R}_+)\right)$.
\end{example}

Introduce now the $C^*$-subalgebra
$\mathcal{T}\left(T_{a_+}^{(0)}\right)$ of the algebra
$\mathcal{T}_0\left(L_\infty^{\{0,+\infty\}}(\mathbb{R}_+)\right)$
which is generated by identity and the Toeplitz operator
$T_{a_+}^{(0)}$ with symbol $a_+=\chi_{[0,1/2]}$. By
Example~\ref{exm12} the corresponding function $\gamma_{a_+,0}\in
L_\infty^{\{0,+\infty\}}(\mathbb{R}_+)$ is continuous on
$[0,+\infty]$, therefore $T_{a_+}^{(0)}$ is self-adjoint and
$\textrm{sp}\, T_{a_+}^{(0)} = \textrm{Range}\,\gamma_{a_+,0} =
[0,1]$. Also, the function $\gamma_{a_+,0}$ is strictly increasing
real-valued function which separates the point of $[0,+\infty]$.
Thus, the algebra $\mathcal{T}\left(T_{a_+}^{(0)}\right)$ consists
of all operators of the form $h\left(T_{a_+}^{(0)}\right)$ with
$h\in C[0,1]$ by functional calculus for $C^*$-algebras.
Summarizing, we have

\begin{theorem}
The algebra
$\mathcal{T}_0\left(L_\infty^{\{0,+\infty\}}(\mathbb{R}_+)\right)
= \mathcal{T}\left(T_{a_+}^{(0)}\right)$, and
\begin{itemize}
\item[(i)] is isomorphic and isometric to $C[0,+\infty]$;
\item[(ii)] is generated by identity and the single Toeplitz
operator $T_{a_+}^{(0)}$; \item[(iii)] consists of all operators
of the form $h\left(T_{a_+}^{(0)}\right)$ with $h\in C[0,1]$.
\end{itemize}
\end{theorem}

We have chosen Toeplitz operator $T_{a_+}^{(0)}$ as the starting
operator because in this specific case the equation
$x=\gamma_{a_+,0}(\xi) = 1-\mathrm{e}^{-\xi}$ admits an explicit
solution. But we can start from any operator
$T_{\chi_{[0,\lambda]}}^{(0)}$ with symbol
$a(u)=\chi_{[0,\lambda]}(u)$, $\lambda\in\mathbb{R}_+$. Indeed,
the function $\gamma_{\chi_{[0,\lambda]},0}(\xi)$ is strictly
increasing which implies that the function
$$\Delta_\lambda(x)=1-(1-x)^{2\lambda}, \quad x\in [0,1],$$ is
strictly increasing as well, and thus the function
$\Delta_\lambda^{-1}$ is well defined and continuous on $[0,1]$.
Clearly, for $\lambda_1, \lambda_2\in\mathbb{R}_+$ we have
$$\left(\Delta_{\lambda_2}\circ
\Delta_{\lambda_1}^{-1}\right)\left(T^{(0)}_{\chi_{[0,\lambda_1]}}\right)
= T^{(0)}_{\chi_{[0,\lambda_2]}},$$ where $\circ$ is the usual
composition of real functions. This means that for any symbol
$a=a(u)=\chi_{[0,\lambda]}(u)$ the operator $T_a^{(0)}$ belong to
the algebra $\mathcal{T}\left(T_{a_+}^{(0)}\right)$, and is the
function of $T_{a_+}^{(0)}$, i.e.,
$T_{\chi_{[0,\lambda]}}^{(0)}=\Delta_\lambda\left(T_{a_+}^{(0)}\right)$.
In fact, the (localization) operator $T^{(0)}_{a_+}$ gives a
reconstruction of a signal on the segment
$\Omega_{1/2}=\mathbb{R}\times (0,1/2]$ and level $0$. Then we may
obtain any operator $T_{\chi_{[0,\lambda]}}^{(0)}$ (giving a
reconstruction of a signal on the segment
$\Omega_\lambda=\mathbb{R}\times (0,\lambda]$ and level $0$) from
this operator $T_{a_+}^{(0)}$, i.e., \textsf{from the
reconstruction of a signal on the segment $\Omega_{1/2}$ and level
$0$ we may obtain a reconstruction of the same signal on the same
level on an arbitrary segment $\Omega_\lambda$} using the function
$\Delta_\lambda$ which is easy to compute.

In particular and quite surprisingly, each Toeplitz operator
$T_a^{(k)}$ with symbols $a\in
L_\infty^{\{0,+\infty\}}(\mathbb{R}_+)$ is a certain continuous
function of the initial operator and this function can be figured
out.

\begin{theorem}\label{thmexchange2}
For each $a=a(u)\in L_\infty^{\{0,+\infty\}}(\mathbb{R}_+)$ the
Toeplitz operator $T_a^{(k)}$ belongs to the algebra
$\mathcal{T}\left(T_{a_+}^{(0)}\right)$, and is the following
function of the operator $T_{a_+}^{(0)}$
$$\left(\nabla_{a,\lambda}^{(k)}\circ\Delta_\lambda\right)\left(T_{a_+}^{(0)}\right)
= T_{a}^{(k)},$$ where
$$\nabla_{a,\lambda}^{(k)}(x)=-\frac{1}{\lambda}\ln(1-x)\int_{\mathbb{R}_+} a(u)(1-x)^{u/\lambda}
L_k^2\left(-\frac{u}{\lambda}\ln(1-x)\right)\,\mathrm{d}u$$ with
$\lambda\in\mathbb{R}_+$ and $x\in [0,1]$.
\end{theorem}

Theorem~\ref{thmexchange2} states that if we know the
reconstruction of a signal on a segment $\Omega_\lambda$ and level
$0$, we might get an arbitrary reconstruction of the signal (as
its filtered version using a real bounded function $a$ of scale
having limits in critical points of boundary of $\mathbb{R}_+$
such that the corresponding function $\gamma_\cdot$ separates the
points of $\overline{\mathbb{R}}_+$) on an arbitrary level $k$
using the function $\nabla_{a,\lambda}^{(k)}$. Theoretically, for
the purpose to study localization of a signal in the time-scale
plane the result of Theorem~\ref{thmexchange2} suggests to
consider certain "nice" symbols on the first level $0$ (indeed,
Toeplitz operators on $\mathcal{A}_2(\Pi)$ with symbols as
characteristic functions of some interval in $\mathbb{R}_+$)
instead of possibly complicated
$L_\infty^{\{0,+\infty\}}(\mathbb{R}_+)$-symbols with respect to
"different microscope" represented by the level $k$. On the other
hand, to compute the corresponding function
$\nabla_{a,\lambda}^{(k)}$ need not be always easy.

\subsection{Unitarily equivalent images of Toeplitz operators for
symbols depending on $\Re\zeta$}\label{sectionFredholm}

For general symbols $a=a(u,v)$ the operator $T_a^{(k)}$ is no
longer unitarily equivalent to a multiplication operator
$\mathfrak{A}_a^{(k)}$. For symbols depending on the second
(vertical) variable in the upper half-plane $\Pi$ a certain class
of pseudo-differential operaors appears. In what follows
$\mathbb{R}_+^2:=\mathbb{R}_+\times\mathbb{R}_+$.

\begin{theorem}[\cite{hutnikIII}, Theorem 3.3]\label{CTO2} Let $(u,v)\in \mathbb{G}$.
If a measurable function $b=b(v)$ does not depend on $u$, then
$T_{b}^{(k)}$ acting on $A^{(k)}$ is unitarily equivalent to the
operator $\mathfrak{B}_b^{(k)}$ acting on $L_2(\mathbb{R}_+)$
given by
$$\left[\mathfrak{B}_b^{(k)}f\right](\xi) = \int_{\mathbb{R}_+}
B_k(\xi,t)\, \hat{b}(\xi-t) f(t)\,\mathrm{d}t, \quad
\xi\in\mathbb{R}_+,
$$ where the function $B_k: \mathbb{R}_+^2\to\mathbb{C}$ has the form
\begin{equation}\label{function B_k}
B_k(\xi,t) = \frac{2\sqrt{t\xi}}{t+\xi}\,
P_k\left(\frac{8t\xi}{(t+\xi)^2}-1\right)\end{equation} with
$$P_n(x) = \frac{1}{2^n n!}\, \frac{\mathrm{d}^n}{\mathrm{d}x^n}
(x^2-1)^n$$ being the Legendre polynomial of degree
$n\in\mathbb{Z}_+$ for $x\in [-1,1]$.
\end{theorem}

Immediately, for each $k\in\mathbb{Z}_+$ the function $B_k:
\mathbb{R}_+^2 \to \mathbb{R}$ has the following remarkable
properties:
\begin{itemize}
\item[(i)] $B_k$ is continuous and bounded on $\mathbb{R}_+^2$;
\item[(ii)] $B_k$ is a symmetric function, i.e.,
$B_k(\xi,t)=B_k(t,\xi)$ for each $(\xi,t)\in\mathbb{R}_+^2$;
\item[(iii)] $B_k(\xi,t)\in C^\infty(\mathbb{R}_+^2)$; \item[(iv)]
$B_k$ is homogeneous (of order 0), i.e., for each $\alpha>0$ holds
$B_k(\alpha \xi,\alpha t) = B_k(\xi,t)$ for each $(\xi,
t)\in\mathbb{R}_+^2$; \item[(v)] $B_k(\xi,\xi)=1$ for all
$\xi\in\mathbb{R}_+$.
\end{itemize}

Again, the above result includes the well-known result for
classical Toeplitz operators on the Bergman space
$\mathcal{A}_2(\Pi)$ as a special case. In fact, for $k=0$
Toeplitz operator $T_b$ with a symbol $b=b(\Re \zeta)=b(v)$ acting
on the Bergman space $\mathcal{A}_2(\Pi)$ is unitarily equivalent
to the following integral operator
$$\left[\mathfrak{B}_b f\right](x) = \int_{\mathbb{R}_+}
\frac{2\sqrt{xy}}{x+y}\, K(y-x) f(y)\,\mathrm{d}y, \quad
x\in\mathbb{R}_+,$$ where $K$ is the Fourier transform of the
function $b(-v)$. As we have already mentioned, our case of
Toeplitz operators $T_a^{(k)}$ depending on a "discrete" weight
parameter $k\in\mathbb{Z}_+$ is different from the case of
Toeplitz operators $T_a^{(\lambda)}$ depending on a "continuous"
weight parameter $\lambda\in (-1,+\infty)$ studied
in~\cite{vasilevskibook} for weighted Bergman spaces.

In what follows we deal with the Fredholm theory for Toeplitz
operator algebras on poly-analytic spaces. Thus, consider the
class of integral operators on $\mathbb{R}_+$ of the form
\begin{equation}\label{operator} (Hf)(x)= \int_{\mathbb{R}_+}
h(x,y)K(x-y) f(y)\,\mathrm{d}y, \quad
x\in\mathbb{R}_+,\end{equation} where
\begin{itemize}
\item[(H1)] for the Fourier transform $\hat{K}$ of the function
$K$ holds $|\hat{K}^{(j)}(\omega)|\leq
\frac{C_j}{(1+\omega^2)^{j/2}}$ for each $j\in\mathbb{Z}_+$;
\item[(H2)] there exists limits $\hat{K}_\pm =
\lim\limits_{\omega\to\pm\infty} \hat{K}(\omega)$, and the Fourier
transform of a function $K_0\in L_1(\mathbb{R})$ may be written in
the form
$$\hat{K}_0(\omega) =
\hat{K}(\omega)-\hat{K}_+\chi_+(\omega)-\hat{K}_-\chi_-(\omega);$$
\item[(H3)] $h(x,y)\in C^\infty(\mathbb{R}_+^2)$, and for all
$\alpha>0$ holds $h(\alpha x,\alpha y)=h(x,y)$; \item[(H4)]
$$\int_{\mathbb{R}_+} \frac{|h(1,t)|}{1+t}\,\frac{\mathrm{d}t}{\sqrt{t}}<+\infty.$$
\end{itemize}\noindent With the operator $H$ we associate the function $q_H$ on $\mathbb{R}$ as follows
$$q_H(\lambda) = 
\frac{1}{\mathrm{i}\pi}\int_{\mathbb{R}_+}
t^{-\mathrm{i}\lambda}\,\frac{h(1,t)}{1-t}\,\frac{\mathrm{d}t}{\sqrt{t}},
\quad
\lambda\in\mathbb{R}.$$ 
Denote by $SV^\infty(\mathbb{R}_+)$ the class of functions
$f(x)\in C_b^\infty(\mathbb{R}_+)$ which are slowly varying at
infinity (in the additive sense) and slowly varying at zero (in
the multiplicative sense), see Remark~\ref{remSV}, i.e.,
$\lim_{x\to+\infty} f'(x) =0$ and $\lim_{x\to 0} x f'(x)=0$. Let
$SV(\mathbb{R}_+)$ be the closure of $SV^\infty(\mathbb{R}_+)$ in
$C_b(\mathbb{R}_+)$. Further, denote by $\mathcal{H}$ the
$C^*$-algebra generated by the integral operators $H$ of the
form~(\ref{operator}) with the property
\begin{equation}\label{supq_H}
\sup_{\lambda\in\mathbb{R}}|q_H(\lambda)|<+\infty,
\end{equation}
and the multiplication operators $f(x)I$ with $f\in
SV^\infty(\mathbb{R}_+)$. In such a case the commutator $[f(x)I,
H] = f(x)H-H f(x)I$ is compact on $L_2(\mathbb{R}_+)$.

For each $k\in\mathbb{Z}_+$ the operator $\mathfrak{B}_b^{(k)}$
with $b(v)\in L_\infty(\mathbb{R})$ is of the
form~(\ref{operator}) with $\hat{K}(\omega)=b(-\omega)$ and is
clearly bounded on $L_2(\mathbb{R}_+)$. The associated function
$q_{\mathfrak{B}_b^{(k)}}$ has the form
$$q_{\mathfrak{B}_b^{(k)}}(\lambda) = \frac{2}{\mathrm{i}\pi} \int_{\mathbb{R}_+}
\frac{t^{-\mathrm{i}\lambda}}{(1+t)(1-t)}
P_k\left(\frac{8t}{(t+1)^2}-1\right)\,\mathrm{d}t, \quad
\lambda\in\mathbb{R},$$ and it is possible to prove,
see~\cite{HH}, that for each $k\in\mathbb{Z}_+$
$$\sup_{\lambda\in\mathbb{R}}
\left|q_{\mathfrak{B}_b^{(k)}}(\lambda)\right|<+\infty,$$ which
implies that $\mathfrak{B}_b^{(k)}\in\mathcal{H}$ for each
$k\in\mathbb{Z}_+$ whenever $b=b(v)\in C(\overline{\mathbb{R}})$.
According to Theorem~\ref{thmder} $\gamma_{a,k}(x)\in
SV(\mathbb{R}_+)$ whenever $a=a(u)\in C_b^1(\mathbb{R}_+)$ such
that $\lim\limits_{u\to+\infty} u a'(u) = 0$, and therefore
$\gamma_{a,k}(x)I\in\mathcal{H}$. Thus, denote by
$\widetilde{SV}(\mathbb{R}_+)$ the $C^*$-algebra generated by
functions $a=a(u)\in C_b^1(\mathbb{R}_+)$ with
$\lim\limits_{u\to+\infty} u a'(u) = 0$. For each
$k\in\mathbb{Z}_+$ consider the $C^*$-algebra
$$\mathcal{T}_k=\mathcal{T}_k\left(C(\overline{\mathbb{R}}),
\widetilde{SV}(\mathbb{R}_+)\right)$$ generated by all Toeplitz
operators $T_b^{(k)}$ with $b=b(v)\in C(\overline{\mathbb{R}})$,
and $T_a^{(k)}$ with $a=a(u)\in\widetilde{SV}(\mathbb{R}_+)$
acting on $A^{(k)}$, where $\zeta=(u,v)\in \mathbb{G}$. It is
immediate that $\mathcal{T}_k$ provides a parameterized family of
operator algebras with compact commutator property and non-compact
semi-commutator property.

Now, introduce two ideals of the algebra $SV(\mathbb{R}_+)$,
\begin{align*}
C_0^{0}(\mathbb{R}_+) & = \Bigl\{a(u)\in SV(\mathbb{R}_+);
\,\,\,\,\lim_{u\to 0} a(u)=0\Bigr\}; \\
C_0^{+\infty}(\mathbb{R}_+) & = \Bigl\{a(u)\in SV(\mathbb{R}_+);
\,\,\,\,\lim_{u\to+\infty} a(u)=0\Bigr\},
\end{align*}and let $\textrm{Sym}\,\mathcal{H}=\mathcal{H}/\mathcal{K}=\widehat{\mathcal{H}}$ be the
Fredholm symbol algebra of the algebra $\mathcal{H}$, where
$\mathcal{K}$ is the ideal of all compact operators on
$L_2(\mathbb{R}_+)$. For each $\xi_0\in\overline{\mathbb{R}}_+$
the local algebra is defined as
$\mathcal{H}(\xi_0)=\widehat{\mathcal{H}}/J(\xi_0)$, where
$J(\xi_0)$ is the closed two-sided ideal of the algebra
$\widehat{\mathcal{H}}$ generated by the maximal ideal of
$C(\overline{\mathbb{R}}_+)$ corresponding to the point $\xi_0$.
Clearly, $C(\overline{\mathbb{R}}_+)$ is a central commutative
subalgebra of $\widehat{\mathcal{H}}$. Denote by $\mathfrak{S}$
the $C^*$-algebra of all vector-valued functions $\sigma$
continuous on $\overline{\mathbb{R}}_+$, where
$\sigma(\xi)\in\mathcal{H}(\xi)$ for each
$\xi\in\overline{\mathbb{R}}_+$, with point-wise operations, and
the norm $\|\sigma\| = \sup_{\xi\in\overline{\mathbb{R}}_+}
\|\sigma(\xi)\|$.

\begin{theorem}
For each $k\in\mathbb{Z}_+$ the Fredholm symbol algebra
$\mathrm{Sym}\,\mathcal{T}_k = \mathcal{T}_k/\mathcal{K}$ of
Toeplitz operator algebra $\mathcal{T}_k$ is isomorphic and
isometric to the algebra $\mathfrak{S}$. The symbol homomorphism
$$\mathrm{sym}_k: \mathcal{T}_k \longrightarrow
\mathrm{Sym}\,\mathcal{T}_k=\mathfrak{S}$$ is generated by the
following mapping of the generators of the algebra $\mathcal{H}$
\begin{align*}
\mathrm{sym}_k: T_a^{(k)} & \longmapsto
\begin{cases}\gamma_{a,k}(\xi)+C_0^0(\mathbb{R}_+), & \xi=0 \\
(\gamma_{a,k}(\xi), \gamma_{a, k}(\xi)), &
\xi\in\mathbb{R}_+; \\
\gamma_{a,k}(\xi)+C_0^{+\infty}(\mathbb{R}_+), &
\xi=+\infty\end{cases}
\\
\mathrm{sym}_k: T_b^{(k)} & \longmapsto
\begin{cases}\frac{1}{2}\left[\Bigl(b(-\infty)+b(+\infty)\Bigr)+
\Bigl(b(-\infty)-b(+\infty)\Bigr)q_{\mathfrak{B}_b^{(k)}}(\lambda)\right], & \xi=0 \\
(b(-\infty), b(+\infty)), &
\xi\in\mathbb{R}_+; \\
b(-\omega), & \xi=+\infty\end{cases}
\end{align*}where $a=a(u)\in\widetilde{SV}(\mathbb{R}_+)$ and
$b=b(v)\in C(\overline{\mathbb{R}})$ with $\zeta=(u,v)\in
\mathbb{G}$.
\end{theorem}

\section{Conclusion}

In this paper we review our recent results on Toeplitz operators
on wavelet subspaces with respect to a special parameterized
family of wavelets from Laguerre functions. As it was shown
in~\cite{abreu} such spaces are in fact the true-poly-analytic
Bergman spaces over the upper half-plane $\Pi$ providing thus a
useful and interesting tool for investigating basic properties of
Toeplitz operators acting on them and their algebras.

We list here some questions and possible directions which could be
interesting for further study.

(i) We suppose that considering an appropriate family of
admissible wavelets $\psi^{(k,\alpha)}$ related to generalized
Laguerre functions $$\ell_k^{(\alpha)}(x) :=
\left[\frac{k!}{\Gamma(k+\alpha+1)}\right]^{1/2}
x^{\alpha/2}\,\mathrm{e}^{-x/2}L_k^{(\alpha)}(x), \quad x\in
\mathbb{R}_+,$$ we may obtain analogous structural results and
representations of weighted (true) poly-analytic Bergman spaces.
\textit{What are the corresponding Toeplitz (Hankel, etc.)
operators acting on these spaces and which properties do they
have}?

(ii) In connection with the results of this paper we may ask
\textit{what happens to properties of Toeplitz operators}
$T_a^{(k)}$ acting on true-poly-analytic Bergman spaces (wavelet
subspaces $A^{(k)}$) \textit{when the weight parameter
$k\in\mathbb{Z}_+$ varies}?

(iii) In particular, \textit{to study the spectral properties of a
Toeplitz operator $T_a^{(k)}$} and the related asymptotic
properties of function $\gamma_{a,k}$ \textit{in dependence on}
$k$, and compare their limit behavior under $k\to+\infty$ with
corresponding properties of the initial symbol $a$. Also, the
similar questions may be stated for the above mentioned Toeplitz
operator acting on weighted poly-analytic spaces, but here the
question of a varying parameter is not clear, because in this case
at least two weighted parameters will appear.

(iv) Recently, an interesting result has been obtained
in~\cite{HMV} for the classical Bergman space of analytic
functions on the upper half-plane. As it is already known, the
$C^*$-algebra generated by Toeplitz operators with bounded symbols
(depending on vertical coordinate $u=\mathrm{Im}(\zeta)$ only, the
so-called vertical Toeplitz operators in terminology
of~\cite{HMV}) is isometrically isomorphic to the $C^*$-algebra
generated by the set $\Gamma_0:=\{\gamma_{a,0};\, a\in
L_\infty(\mathbb{R}_+)\}$. In~\cite{HMV} authors showed that
$\Gamma_0$ is dense in the space $\mathrm{VSO}(\mathbb{R}_+)$ of
all very slowly oscillating functions on the positive half-line.
Thus, we conjecture the following more general result: \textit{for
each $k\in\mathbb{Z}_+$ the set}
$$\Gamma_k:=\{\gamma_{a,k};\, a\in L_\infty(\mathbb{R}_+)\}$$
\textit{is dense in $\mathrm{VSO}(\mathbb{R}_+)$.}

(v) In practical applications certain algebraic operations with
symbols and operators naturally appear. In signal analysis, the
problem of finding a filter that has the same effect as two
filters arranged in series amounts to the computation of the
product of two localization operators. Thus, \textit{what is the
product of two Toeplitz operators (in exact, or at least
approximate formulas)}? The answer does not seem to be so simple
and seems to depend on the availability of a useful formula for
Toeplitz operator. We think that some technical information in
this direction obtained in this paper studying the particular
cases of symbols on the affine group $\mathbb{G}$ (e.g., if $a$
depends only on horizontal variable $u\in\mathbb{R}_+$ in the
upper half-plane, then $T_a^{(k)}$ is a Fourier multiplier) may be
helpful and crucial.

(vi) The continuous wavelet transform in the one-dimensional case
can be obtained in two ways: one from the theory of
square-integrable group representation, and the other from the
Calder\'on representation formula. Also, it is known that in the
one-dimensional case these two different ways can induce the same
results. However, in the higher dimensional case these two ways
will induce two different results. One is the Calder\'on
representation formula, which induces a decomposition of
$L_2(\mathbb{R}_+\times \mathbb{R}^n,
u^{-n-1}\mathrm{d}u\mathrm{d}v)$, and the other is the wavelet
transform associated with the square-integrable group
representation. Also, there exist other ways how to extend wavelet
analysis to higher dimensions, cf.~\cite{AAGbook}. Each such a
case generates its own (possibly different) class of localization
operators. \textit{What is the "natural" extension of our results
for Toeplitz operators to higher dimensions}?

Immediately, there are many other questions dealing with various
contexts, e.g., in quantization problems, (discrete and
continuous) frame theory, engineering applications, etc. We hope
this paper will stimulate a further interest and development in
this topic of intersection of poly-analytic function theory and
time-scale (or, more generally, time-frequency) analysis.

\section*{Acknowledgement}

This paper was developed as a part of the project named "Centre of
Excellence for Integrated Research \& Exploitation of Advanced
Materials and Technologies in Automotive Electronics", ITMS
26220120055. The first author acknowledges a partial support of
grant VEGA 2/0090/13.


\vspace{5mm}

\noindent \small{Ondrej Hutn\'ik, Institute of Mathematics,
Faculty of Science, Pavol Jozef \v Saf\'arik University in Ko\v
sice, {\it Current address:} Jesenn\'a 5, SK 040~01 Ko\v sice,
Slovakia,
\newline {\it E-mail address:} ondrej.hutnik@upjs.sk}

\vspace{5mm}

\noindent \small{M\'aria Hutn\'ikov\'a, Department of Physics,
Faculty of Electrical Engineering and Informatics, Technical
University of Ko\v sice, {\it Current address:} Park Komensk\'eho
2, SK 042~00 Ko\v sice, Slovakia,
\newline {\it E-mail address:} maria.hutnikova@tuke.sk}

\end{document}